\documentclass[11pt]{article}
\usepackage[greek,english]{babel}
\usepackage[iso-8859-7]{inputenc}
\usepackage{amssymb}
\usepackage{amsmath}
\usepackage{amsthm}
\usepackage{latexsym}
\usepackage{amsfonts}
\usepackage{graphicx}
\usepackage{graphics}

\theoremstyle{plain}
\newtheorem{thm}{Theorem}[section]   
\newtheorem{prop}[thm]{Proposition}
\newtheorem{lem}[thm]{Lemma}

\theoremstyle{definition}

\newtheorem*{Proof}{Proof}

\newcommand{\bbb}[1]{\mbox{\boldmath$#1$}}

\newcommand{\OO} {{\varOmega}}

\newcommand{\ga} {{\gamma}}

\newcommand{\Ga} {{\varGamma}}

\newcommand{\ld} {{\ldots}}
\newcommand{\sm} {{\smallsetminus}}
\newcommand{\thi} {{\theta}}

\newcommand{\de} {{\delta}}
\newcommand{\De} {{\varDelta}}

\newcommand{\la} {{\lambda}}
\newcommand{\el} {{\ell}}

\newcommand{\Fi} {{\varPhi}}
\newcommand{\e} {{\varepsilon}}

\newcommand{\dis}{\displaystyle}
\newcommand{\ssum}{\sum\limits}

\newcommand{\cp}{{\cal{P}}}
\newcommand{\cl}{{\cal{L}}}

\newcommand{\tg}{{\widetilde{g}}}

\newcommand{\mfd}{{\mathfrak{D}}}
\newcommand{\ra}{{\rightarrow}}
\newcommand{\oD}{{\overline{D}}}
\newcommand{\oV}{{\overline{V}}}
\newcommand{\oOO}{{\overline{\varOmega}}}

\newcommand{\lo}{{\overset{\circ}{L}}}

\newcommand{\tint}{{\text{Int}}}
\newcommand{\tind}{{\text{Ind}}}
\newcommand{\qb}{$\quad\blacksquare$}
\def\1{\it1\hspace*{-0.150cm}{\footnotesize{I}}}

\def\R{{\mathbb{R}}}
\def\C{{\mathbb{C}}}

\def\N{{\mathbb{N}}}

\begin{document}
\title{\bf A general lower bound for the asymptotic convergence factor}
\author{\bf N. Tsirivas}\footnotetext{\hspace{-0.5cm}}
\footnotetext{{The research project is implemented within the framework of the Action ``Supporting Postdoctoral Researchers'' of the Operational Program ``Educational and Lifelong Learning'' (Action's Beneficiary: General Secretariat for Research and Technology), and is co-financed by the European Social Fund (ESF) and the Greek State.}}
\date{}
\maketitle
\noindent
{\bf Abstract:} We provide a rather general and very simple to compute lower bound for the asymptotic convergence factor of compact subsets of $\C$ with connected complement and finitely many connected components. 

%
\setcounter{section}{-1}
\section{Introduction}\label{sec0}
\noindent

The subject of this work has many connections with the theory of approximation \cite{4}, the problem of solving large-scale matrix problems by Krylov subspace iterations, and dijital filtering.

Our point of view is a classical problem of approximation.

Let us see how this problem is arises.

Suppose that $L$ is a non-connected compact subset of $\C$ with connected complement.

Hence, $L:=\bigcup\limits^m_{i=0}K_i$, for some $m\in\N$, $m\ge1$ where $K_i$, $i=0,1,\ld,m$ are the connected components of $L$. We also assume that each component $K_i$ does not reduce to a point. Let $p_i$, $i=0,1,\ld,m$ be $m+1$ different complex polynomials, that is $p_i\neq p_j$ for every $i,j\in\{0,1,\ld,m\}$, $i\neq j$.

Consider the function $F:L\ra\C$, that is defined by the formula:
\[
F(z)=p_j(z) \ \ \text{if} \ \ z\in K_j, \ \ \text{for every} \ \ j\in\{0,1,\ld,m\}.
\]
We fix some positive number $\de$. The problem is to find explicitly a polynomial $p$ such that $\|F-p\|<\de$ and also to find the relation between $p$ and $\de$, if possible. It turns out that the notion of asymptotic convergence factor for a compact set is extremely useful in studying the above problem. Based on the previous notation we proceed with the relevant definition. 

For every $n=1,2,\ld$, let $V_n$ be the set of complex polynomials with degree at most $n$.

We denote
\[
d_n, F:=\min\{\|F-p\|,\;p\in V_n\} \ \ \text{for} \ \ n=1,2,\ld.
\]
Of course for every $n\in\N$, there exists some $p\in V_n$ such that $d_n,F=\|F-p\|_L$, and the polynomial $p$ is unique (\cite{14}) for every $n\ge1$. Even, if the formulation of the problem of finding the above best polynomial $p$ that minimizes the quantity $\|F-p\|_L$ is simple this is usually unknown and difficult to compute (see \cite{6}, page 11).

However, if the compact set $L$ has a simple structure or good regularity properties the previous approximation problem can be solved. Despite this, the computation of the best polynomial is difficult even in simple cases, for instance the union of two disjoint closed disks, and in most of the cases this is done with complicated numerical methods \cite{6}.

A classical theorem in this area, see \cite{14}, \cite{6}, is the following.
\begin{thm}\label{thm1.1}
The number $\rho_L:=\underset{n\ra+\infty}{\lim\sup}\,d^{\frac{1}{n}}_n,F$ is a positive constant such that $\rho_L\in(0,1)$, it is independent from the function $F$ and it is dependent only on the compact set $L$.
\end{thm}

The number $\rho_L$ is called the asymptotic convergence factor of $L$ and is a characteristic for the compact set $L$. Of course the knowledge of the above number $\rho_L$ for $L$ is a crucial point for the solution of the initial approximation problem.

However, the number $\rho_L$ is very difficult to be computed in general, \cite{6}, \cite{10}, \cite{11}. So, it is desirable for a simple compact set $L$ to obtain ``good" estimates from above and below of the number $\rho_L$.

In this paper we give an easily computed lower bound for the number $\rho_L$, which is best possible in a certain sense. Our main result is the following

\begin{thm} \label{mainthm}
Under the above assumptions and notation we have 
$$\rho_L \geq \max_{j=0,\ldots ,m}\sup_{z\in K_j^o}\frac{ \textrm{dist}(z,K_j^c)}{\textrm{dist}(z,L\setminus K_j)} .$$
\end{thm}

In order to prove this theorem we introduce an other characteristic number $\thi_L$ for a compact set $L$, which is defined in a bit complicated way in the next section. However, this number gives us the necessary potential theoretic tools in order to establish that $\rho_L=\thi_L$, that is $\thi_L$ is nothing but the asymptotic convergence factor. The potential theoretic view, assigns to the number $\rho_L$ new important properties. We also provide a variety of simple examples of certain compact sets where the number $\thi_L$ can be computed by a very simple algebraic formula and not with numerical methods.

Results concerning the computation of $\rho_L$ for certain compact set $L$ can be found in \cite{1}, \cite{6}, \cite{10}, \cite{11}.\vspace*{0.2cm} \\
\noindent
{\bf Remark.} We note that Proposition \ref{prop2.3} of this paper is used in a substantial way in order to prove the main result Theorem \ref{thm1.1} of \cite{13}. The theme of \cite{13} is related to universal Taylor series. On the other hand, according to Remark 3.3 of \cite{13}, half of the main result of this paper, namely Theorem \ref{thm1.4}, can be deduced by a completely different method based entirely on results of universality for Taylor series \cite{13}. This means that there exists a close relation between the results of this paper and \cite{13}. 
\section{The number $\thi_L$ and its lower bound }\label{sec1}
\noindent

We begin with the necessary terminology. For the topological concepts of this paper we refer to the classical book of Burckel \cite{3}.

More specific for the definitions of a curve, or a loop, or an arc, or a simple curve, or a smooth curve see Definition 1.11 \cite{3}.

For the definition of a simply-connected subset $A$ of $\C$ see Definition 1.36 \cite{3}. With a Jordan curve we mean a homeomorphism in $\C$ of a circle. If $\ga$ is a smooth Jordan curve and $w\in\C\setminus\ga$, the index $\tind_\ga(w):=\dfrac{1}{2\pi i}\dis\int_\ga\dfrac{1}{z-w}dz$.

For a compact subset $K$ of $\C$ and a Jordan curve $\ga$ we write
\[
\tind_\ga(K):=\{\tind_\ga(w),\;w\in K\},
\]
when $\ga\cap K=\emptyset$.

The definition of interior, $\tint(\ga)$ and Exterior $Ex(\ga)$ of a Jordan curve $\ga$ is given in Definition 4.45 (i) of \cite{3}. For results about potential theory we refer to the classical books \cite{2} and \cite{8}.

Below we prove a series of useful topological lemmas and the main result of this paper, that is a simple estimation of the lower bound of the number $\rho_L$, for many cases of compact sets.
\begin{lem}\label{lem1.1}
Let $V\subseteq\C$, $V\neq\C$ be a simply connected domain, $K\subseteq V$, $K$ compact, set.

Then, there exists a smooth Jordan curve $\ga\subset V$ such that $\tind_\ga(K)=\{1\}$.
\end{lem}
\begin{Proof}
Let $D$ be the open unit disc. By the Riemann mapping theorem there exists a conformal mapping $f:D\ra V$, that is $1-1$ and onto. We set $L:=f^{-1}(K)$. Of course the set $L$ is a compact subset of $D$. Let $r_0\in(0,1)$ such that $L\subset D(0,r_0)$, where
\[
D(0,r_0):=\{z\in\C\mid\,\mid z\mid<r_0\}.
\]
We set
\[
\Ga:=C(0,r_0):=\{z\in\C\mid\,\mid z\mid=r_0\}.
\]
We consider the circle $\ga_0:[0,1]\ra\C$, $\ga_0(t)=r_0e^{2\pi it}$, $t\in[0,1]$ and we set $\ga:=f\circ\ga_0$. We set $\ga^\ast:=\ga([0,1])$ and we write simply $\ga^\ast=\ga$ without confusion.

It is easy to show that the curve $\ga$ is a smooth Jordan curve such that $K\cap\ga=\emptyset$. So the number $\tind_\ga(w)$ has sense for every $w\in K$.

We fix some $w_0\in K$. We compute the number $\tind_\ga(w_0)$.

We have:
\begin{align}
\tind_\ga(w_0):&=\frac{1}{2\pi i}\int_\ga\frac{1}{z-w_0}dz \nonumber\\
&=\frac{1}{2\pi i}\int^1_0\frac{1}{\ga(t)-w_0}\cdot\ga'(t)dt \nonumber\\
&=\frac{1}{2\pi i}\int^1_0\frac{f'(\ga_0(t))\cdot\ga'_0(t)}{f(\ga_0(t))-w_0}dt.  \label{eq1}
\end{align}
We consider the function $g:D\setminus\{z_0\}\ra\C$ defined by the formula:

\[
g(z):=\frac{f'(z)}{f(z)-w_0}, \ \ \text{where} \ \ z_0:=f^{-1}(w_0), \ \ z\in D\setminus\{z_0\}.
\]
Obviously, the function $g$ is well defined and holomorphic in $D\setminus\{z_0\}$ and has a singularity in $z_0$. Obviously $z_0$ is pole of $g$. It holds $z_0\notin\ga_0$, so the integral $\int\limits_{\ga_0} g(z)dz$ is well defined.

Now, we have
\begin{align}
\int_{\ga_0}g(z)dz&=\int^1_0g(\ga_0(t))\cdot\ga'_0(t)dt \nonumber \\
&=\int^1_0\frac{f'(\ga_0(t))\cdot\ga'_0(t)}{f(\ga_0(t))-w_0}.  \label{eq2}
\end{align}
By (\ref{eq1}) and (\ref{eq2}) we take:
\begin{eqnarray}
\tind_\ga(w_0)=\frac{1}{2\pi i}\int_{\ga_0}\frac{f'(z)}{f(z)-w_0}dz.  \label{eq3}
\end{eqnarray}
Because $L\subset D(0,r_0)$ by definition we have that $z_0\in\tint(\ga_0)$.

We compute easily that $\dis\lim_{z\ra z_0}(z-z_0)g(z)=1$. This gives that $z_0$ is a simple pole for $g$ and $Res(g,z_0)=1$. So by the calculus of residues we take
\begin{eqnarray}
\int_{\ga_0}g(z)dz=2\pi i\,Res(g,z_0)\cdot\tind\ga_0(z_0)=2\pi i.  \label{eq4}
\end{eqnarray}
By (\ref{eq3}) and (\ref{eq4}) we take $\tind\ga(w_0)=1$ and the result follows. \qb
\end{Proof}

Using Lemma \ref{lem1.1} we prove now the following lemma, for compact subsets of $\C$, with finite many connected components only.
\begin{lem}\label{lem1.2}
Let $K=\bigcup\limits^m_{i=1}K_i$, be a compact set with connected complement where $K_i$, $i=1,2,\ld,m$ be the connected components of $K$, $m>1$.

Then, there exist smooth Jordan curves $\de_i$, $i=1,2,\ld,m$, pairwise disjoint such that $\tind_{\de_i}(K_i)=\{1\}$ and $K_i\subset\tint(\de_i)$, for $i=1,2,\ld,m$ and every one of them has all the others in its exterior.
\end{lem}
\begin{Proof}
First of all we can choose bounded open subset of $\C$, $G_i$, $i=1,2,\ld,m$, pairwise disjoint such that $K_i\subset G_i$, for $i=1,2,\ld,m$.

By Proposition (iv), page 99 of \cite{3} we have that $K^c_i$ is connected for $i=1,2,\ld,m$. Now, by Costakis and Grosse-Erdmann lemma \cite{5}, \cite{7} we take that there exists open simply connected sets $V_i$, $i=1,2,\ld,m$ such that $K_i\subset V_i\subset G_i$, for $i=1,2,\ld,m$. By Corollary 4.66, page 114 of \cite{3}, we have that if we write $V_i=V^j_i$, $j\in J$ for every $i=1,2,\ld,m$, where $V^j_i$, $j\in J$ to be the connected components of $V_i$, and $J$ is a set of indices then $(V^j_i)^c$, $j\in J$ are connected sets. It is easy to see that there exists unique $j_i\in J$ such that $K_i\subset V^{j_i}_i$ for every $i=1,2,\ld,m$. So, we have that for every $i=1,2,\ld,m$ there exists a bounded and simply connected domain $V^{j_i}_i$ such that $K_i\subset V^{j_i}_i\subset G_i$.

To avoid complicated symbolism we write simply $V_i$ instead of $V^{j_i}_i$. So, we have that for every $i=1,2,\ld,m$ there exists a bounded and simply connected domain $V_i$, such that:
\[
K_i\subset V_i \subset G_i, \qquad i=1,2,\ld,m.
\]
Now, we apply Lemma \ref{lem1.1} and we take that for every $i=1,2,\ld,m$ there exists a smooth Jordan curve $\de_i\subset V_i$ such that: $\tind_{\de_i}(K_i)=\{1\}$, where it is supposed of course that $\de_i\cap K_i=\emptyset$ for $i=1,2,\ld,m$. Because $G_i\cap G_j=\emptyset$ for $i,j\in\{1,2,\ld,m\}$, $i\neq j$ we have that the curves $\de_i$, $i=1,2,\ld,m$ are pairwise disjoint.

Now, let $w\in K_i$, for some $i\in\{1,2,\ld,m\}$. Because $K_i\cap\de_i=\emptyset$ we have that $w_i\in\tint(\de_i)$ or $w_i\in Ex(\de_i)$. If $w_i\in Ex(\de_i)$ then $\tind_{\de_i}(w)=0$, by Theorem 10.10 of \cite{9}. So we have $K_i\subset\tint(\de_i)$ for every $i=1,2,\ld,m$.

Now, let $i,j\in\{1,2,\ld,m\}$, $i\neq j$. We show now that $\de_j\subset Ex(\de_i)$. We have $\de_i\subset V_i$, so $V^c_i\subset\de^c_i$. But $\de^c_i=Ex(\de_i)\cup\tint(\de_i)$. Because $V^c_i$ is connected we have that $V^c_i\subset\tint(\de_i)$ or $V^c_i\subset Ex(\de_i)$. But $V^c_i$ is unbounded, because $V_i$ is bounded and $\tint(\de_i)$ is bounded, so
\setcounter{equation}{0}
\begin{eqnarray}
V^c_i\subset Ex(\de_i).  \label{eq1}
\end{eqnarray}
Now we have
\begin{eqnarray}
V_i\subset G_i\Rightarrow G^c_i\subset V^c_i.  \label{eq2}
\end{eqnarray}
We have also
\begin{eqnarray}
G_i\cap G_j=\emptyset\Rightarrow G_j\subset G^c_i.  \label{eq3}
\end{eqnarray}
Also we have
\begin{eqnarray}
\de_j\subset V_j\subset G_j.  \label{eq4}
\end{eqnarray}
By (\ref{eq1}), (\ref{eq2}), (\ref{eq3}) and (\ref{eq4}) we have: $\de_j\subset Ex(\de_j)$. This gives of course that every one from the smooth Jordan curves $\de_i$, $i=1,2,\ld,m$ has all the others in its exterior and the proof of Lemma \ref{lem1.2} is completed. \qb
\end{Proof}

Now, we fix a compact subset $L$ of $\C$, with connected complement such that
\[
L:=\bigcup^{m_0}_{i=0}K_i, \ \ m_0\in\N, \ \ m_0\ge1, \ \ \text{and} \ \ K_i, \ \ i=0,1,\ld,m_0
\]
to be the connected components of $L$.

We consider the set\\
$\mfd_L:=\{\De\in\cp(\C)\mid$ there exist $m_0+1$ smooth Jordan curves $\de_i$, for $i=0,1,\ld,m_0$ such that $\De=\bigcup\limits^{m_0}_{i=0}\de_i$, $K_i\subset\tint(\de_i)$ and $\tind_{\de_i}(K_i)=\{1\}$ for every $i=0,1,\ld,m_0$ and $\bigcup\limits^{m_0}_{i=0\atop i\neq j}\de_i\subset Ex(\de_j)$ for every $j=0,1,\ld,m_0\}$.

By Lemma \ref{lem1.2} we have $\mfd_L\neq\emptyset$ and by the proof of Lemma \ref{lem1.2} we can show easily that the set $\mfd_L$ is uncountable.

The set $\mfd_L$ is of course well defined and non empty without any other restriction.

From now on we suppose that $\overset{\circ}{K_0}\neq\emptyset$.

Let $\OO:=(\C\setminus L)\cup\{\infty\}$. Then $\OO$ is a proper subdomain of $\C_\infty:=\C\cup\{\infty\}$ and by the fact that $\overset{\circ}{K_0}\neq\emptyset$ we take easily that $\partial\OO$ is non-polar. This gives that there exists the unique Green's function $g_\OO$ for $\OO$, with pole at infinity [Definition 4.4.1 and Theorem 4.4.2,\cite{8}].
Let $\De\in\mfd_L$. We write
\[
\thi_{L,\De}:=\max e^{-g_\OO(z,\infty)}:=\max\{x\in\R\mid\;\exists \;z\in\De:x=e^{-g_\OO(z,\infty)}\}.
\]
It is obvious that the number $\thi_{L,\De}$ is a well defined positive number in $(0,1)$, because $\De$ is a compact set and the green's function $g_\OO$ is continuous in $\OO$ and $\De\subset\OO\setminus\{\infty\}$.

We define
\[
\thi_L:=\inf\{x\in\R\mid\;\exists\;\De\in\mfd_L:x=\thi_{L,\De}\}.
\]
By the above the number $\thi_L$ is a well defined number in $[0,1)$ and corresponds uniquely to the compact set $L$ by its definition.

We fix $z_0\in\overset{\circ}{K_0}$.

We write $L_1:=L\setminus K_0$, $r_0:=dist(z_0,K^c_0)$, $h_0:=dist(z_0,L_1)$.

We have the following lemma.
\begin{lem}\label{lem1.3}
By the above definitions we have easily that:

1) $r_0>0$

2) $r_0<h_0$

3) $h_0<+\infty$.

Let some $\De=\bigcup\limits^{m_0}_{i=0}\de_i\in\mfd_L$.

We fix some $w_0\in K^c_0$ such that: $h_0=|z_0-w_0|$. Obviously there exists the unique $i_0\in\{1,2,\ld,m_0\}$ such that $w_0\in K_{i_0}$. We denote
\[
I:=[z_0,w_0]:=\{z\in\C\mid\;\exists\;t\in[0,1]:z=(1-t)z_0+tw_0\}.
\]
Then we have: $\de_{i_0}\cap I\neq\emptyset$.
\end{lem}
\begin{Proof}
We suppose that $\de_0\cap I=\emptyset$ to take a contradiction. Then
\setcounter{equation}{0}
\begin{eqnarray}
I\subset\de^c_0.  \label{eq1}
\end{eqnarray}
Because $\de_0$ is a Jordan curve we have
\begin{eqnarray}
\de^c_0=\tint(\de_0)\cup Ex(\de_0).  \label{eq2}
\end{eqnarray}
Because the segment $I$ is connected we take by (\ref{eq1}) and (\ref{eq2}) that
\begin{eqnarray}
I\subset\tint(\de_0) \quad \text{or} \label{eq3}
\end{eqnarray}
\begin{eqnarray}
I\subset Ex(\de_0).  \label{eq4}
\end{eqnarray}
We suppose that (\ref{eq4}) holds. Then
\begin{eqnarray}
z_0\in I\subset Ex(\de_0)\Rightarrow z_0\in Ex(\de_0).  \label{eq5}
\end{eqnarray}
We have
\begin{eqnarray}
z_0\in\overset{\circ}{K_0}\subset K_0\subset\tint(\de_0),  \label{eq6}
\end{eqnarray}
by Lemma \ref{lem1.2}.

By (\ref{eq5}) and (\ref{eq6}) we have $\tint(\de_0)\cap Ex(\de_0)\neq\emptyset$ that is false. So we have $I\subset\tint(\de_0)$. Thus
\begin{eqnarray}
w_0\in\tint(\de_0).  \label{eq7}
\end{eqnarray}
Because $w_0\in K_{i_0}\Rightarrow w_0\in K_{i_0}\subset\tint(\de_{i_0})$ by Lemma \ref{lem1.2}, so
\begin{eqnarray}
w_0\in\tint(\de_{i_0}).  \label{eq8}
\end{eqnarray}
By the above we have:
\begin{eqnarray}
\tint(\de_0)\cap\tint(\de_{i_0})\neq\emptyset.  \label{eq9}
\end{eqnarray}
We have $\de_{i_0}\subset Ex(\de_0)$. So
\[
\de_{i_0}\cap\tint(\de_0)=\emptyset\Rightarrow\tint(\de_0)\subset\de^c_{i_0}.
\]
But
\[
\de^c_{i_0}=\tint(\de_{i_0})\cup Ex(\de_{i_0}).
\]
Thus
\begin{eqnarray}
\tint(\de_0)\subset\tint(\de_{i_0}) \ \ \text{or} \ \ \tint(\de_0)\subset Ex(\de_{i_0}).  \label{eq10}
\end{eqnarray}
Because of (\ref{eq9}), relation (\ref{eq10}) gives
\begin{eqnarray}
\tint(\de_0)\subset\tint(\de_{i_0}).  \label{eq11}
\end{eqnarray}
Similarly exactly with (\ref{eq11}) we take:
\begin{eqnarray}
\tint(\de_{i_0})\subset\tint(\de_0).  \label{eq12}
\end{eqnarray}
By (\ref{eq11}) and (\ref{eq12}) we take:
\begin{eqnarray}
\tint(\de_0)=\tint(\de_{i_0}).  \label{eq13}
\end{eqnarray}
By (\ref{eq13}) and Theorem 4.41 of \cite{3} we take:
\begin{eqnarray}
\de_{i_0}=\partial\tint(\de_{i_0})=\partial\tint(\de_0)=\de_0,  \label{eq14}
\end{eqnarray}
that is false because $\de_0\cap\de_{i_0}=\emptyset$. So, we have a contradiction that gives us that $\de_0\cap I\neq\emptyset$. \qb
\end{Proof}

In addition to the above we suppose now that every connected component $K_i$, $i=1,2,\ld,m_0$, of $L$ contains more than one point. We prove now the main result of this paper.
\begin{thm}\label{thm1.4}
With the above notations we have:
\[
\thi_L\ge\frac{r_0}{h_0}.
\]
\end{thm}
\begin{Proof}
We consider the number $w_0\in K^c_0$, $w_0\in K_{i_0}$ for some $i_0\in\{1,\ld,m_0\}$ as in Lemma \ref{lem1.3}. We set
\[
V:=D(z_0,h_0)\setminus K_0=D(z_0,h_0)\cap K^c_0.
\]
The set $V$ is an open subset of $\C$, obviously and it is easy to see that $V\subset\OO$.

We have $w_0\in K^c_0$. So there exists $\e_0>0$ such that $D(w_0,\e_0)\subset K^c_0$. We set
\[
V_1:=D(w_0,\e_0)\cap D(z_0,h_0).
\]
It is easy to see that $\emptyset\neq V_1\subset V$, so $V\neq\emptyset$.

We consider the function $\widetilde{g}:\oOO\ra R$ with the formula:
\[
\tg(z):=\left.\begin{array}{ccc}
                0 & \text{if} & z\in\partial\OO \\
                g_\OO(z,\infty) & \text{if} & z\in\OO
              \end{array}\right\}.
              \]
Using Theorem 4.4.9 of \cite{8} it is easy to see that $\tg$ is continuous. To check this we remark that $\OO$ is a regular domain because the connected components of $L$ are finite and every one from them contains more than one point.

We have that $V\subset\OO$ and the set $V$ is bounded. So the set $\oV$ is a compact subset of $\oOO$ where $\tg$ is continuous, so $\tg$ takes its maximum value on $\oV$ in some point $w_1\in\oV$. So we have:
\[
\tg(w_1)=\max_{z\in\oV}\tg(z).
\]
Because $V\neq\emptyset$ there exists $z_2\in V$, of course $z_2\in\OO$. By Theorem 4.4.3 of \cite{8} we have
\[
g_\OO(z_2,\infty)>0\Rightarrow\tg(z_0)>0.
\]
So $\tg(w_1)\ge\tg(z_2)>0$. But $w_1\in\oV\subset\oOO$. If we had $w_1\in\partial\OO$ then $\tg(w_1)=0$ that is false. So, we have $w_1\in\oOO\setminus\partial\OO\Rightarrow w_1\in\OO$.

We prove now the claim:
\[
\partial V\subset C(z_0,h_0)\cup\partial\OO.
\]
Let some $j_0\in\partial V$. Because $V\subset\OO\Rightarrow\oV\subset\oOO$ and $\partial V\subset\oV\Rightarrow\partial V\subset\oOO$.

We suppose that $j_0\notin\partial\OO$. Because $j_0\in\oOO$ we take $j_0\in\OO$. So we have $j_0\notin K_0$. But we have
\[
V=D(z_0,h_0)\setminus K_0\subset D(0,h_0)\Rightarrow\oV\subset\overline{D(z_0,h_0)}.
\]
If we have $j_0\in D(z_0,h_0)$, then we have $j_0\in D(z_0,h_0)\setminus K_0=V$. But we have $\partial V\cap V=\emptyset$, because $V$ is open, and $j_0\in\partial V\cap V$, that is false. So $j_0\notin D(z_0,h_0)$.

Thus
\[
j_0\in\oV\subset\overline{D(0,h_0)} \ \ \text{and} \ \ j_0\notin D(z_0,h_0).
\]
So
\[
j_0\in C(z_0,h_0)\Rightarrow|j_0-z_0|=h_0,
\]
and the proof of this claim is complete.

We suppose that $w_1\in V$. Let $V_{w_1}$ to be the unique connected component of $V$ such that $w_1\in V_{w_1}$. By the maximum principle Theorem 1.1.8, (a) \cite{8} we have that $\tg$ is a constant on $V_{w_1}$. Because $V_{w_1}$ is an open set let $\e_2>0$ such that $D(w_1,\e_2)\subset V_{w_1}$. Then $\tg(z)=\tg(w_1)$ for every $z\in D(w_1,\e_2)$ of course. By the identity principle, Theorem 1.1.7 \cite{8} we have that $g$ is a constant in $\OO$, that is false because $\dis\lim_{z\ra\infty}g(z)=\infty$.

So, we have $w_1\notin V$, but $w_1\in\oV$. Thus, we have $w_1\in\partial V$. Be the previous claim we have that $w_1\in C(z_0,h_0)\cup\partial\OO$. But $w_1\in\OO\Rightarrow w_1\notin\partial\OO$. So, we have that $w_1\in C(z_0,h_0)\Rightarrow|w_1-z_0|=h_0$ and $w_1\in\OO$. So we have proved that:
\[
g_\OO(w_1,\infty)>0, \ \ g_\OO(w_1,\infty)=\max_{z\in\oV}\tg(z) \ \ \text{and} \ \ w_1\in C(z_0,h_0)\cap\OO.
\]
It is easy to check that $\emptyset\neq\de_0\cap I\subset V$.

We set
\[
C_1:=\{z\in\C\mid\,\mid z-z_0\mid\le r_0\}.
\]
It is easy to see that $C_1\subset K_0\Rightarrow C_1\subset\bigcup\limits^{m_0}_{i=0}K_i=L$. We set $C^c_1=\OO_1$. We take $\OO\subset\OO_1$, where we suppose that $\infty\in\OO_1$. By Corollary 4.4.5 \cite{8} we have
\[
g_\OO(z,\infty)\le g_{\OO_1}(z,\infty), \ \ \text{for every} \ \ z\in\OO.
\]
So we have $g_\OO(w_1,\infty)\le g_{\OO_1}(w_1,\infty)$.

Let some $z_1\in V$. Then
\[
g_\OO(z_1,\infty)\le g_\OO(w_1,\infty).
\]
So
\begin{align*}
g_\OO(z_1,\infty)&\le g_{\OO_1}(w_1,\infty)=\log\bigg(\frac{|w_1-z_0|}{r_0}\bigg)=\log\bigg(\frac{h_0}{r_0}\bigg)
\Rightarrow e^{g_\OO(z_1,\infty)}\le\frac{h_0}{r_0}\Rightarrow\frac{r_0}{h_0}\\
&\le e^{g_\OO(z_1,\infty)}\le\max_{z\in\de_0}e^{g_\OO(z,\infty)}\le\max_{z\in\De}
e^{g_\OO(z,\infty)}=\thi_{L,\De}.
\end{align*}
So we have:
\[
\frac{r_0}{h_0}\le\thi_{L,\De} \ \ \text{for every} \ \ \De\in\cl_L.
\]
Thus, we take $\dfrac{r_0}{h_0}\le\thi_L$ and the proof of this lemma is complete. \qb
\end{Proof}

Theorem \ref{thm1.4} gives us a simple lower bound for the number $\thi_L$.

We will prove that in some cases this lower bound is optimal in some sense.

More specifically:

Let $D$ be the open unit disc and we denote $K_0:=\oD=\{z\in\C\mid\,\mid z\mid\le1\}$ the closed unit disc, for the sequel.

We fix some positive natural number, $h_0>1$. We set:
\[
C_{h_0}:=\{L\subseteq\C\mid L
\]
is compact with connected complement, $L=\bigcup\limits^m_{i=0}K_i$, $m\ge1$, $K_0:=\oD$ where $K_i$, $i=0,1,\ld,m$ to be the connected components of $L$ and $K_i$, $i=1,\ld,m$ contains more than one point and $dist(\{0\},L\setminus K_0)=h_0$, $h_0\in L\}$.

Of course, by Theorem \ref{thm1.4} we have:
\[
\thi_L\ge\frac{1}{h_0} \ \ \text{for every} \ \ L\in C_{h_0}.
\]
We prove the following proposition.
\begin{prop}\label{prop1.5}
It holds $\inf\{\thi_L,\;L\in C_{h_0}\}=\dfrac{1}{h_0}$.
\end{prop}
\begin{Proof}
We set $I:=\inf\{\thi_L,\;L\in C_{h_0}\}$. Of course, we have
\setcounter{equation}{0}
\begin{eqnarray}
I\ge\dfrac{1}{h_0}, \label{eq1}
\end{eqnarray}
by Theorem \ref{thm1.4}.

It holds
\[
1-\frac{1}{h_0}>0.
\]
We prove that for every $\de\in\Big(0,1-\dfrac{1}{h}\Big)$ there exists some $L'\in C_{h_0}$ such that
\[
\thi_{L'}<\de+\frac{1}{h_0},
\]
that yields of course that
\begin{eqnarray}
I\le\frac{1}{h_0},  \label{eq2}
\end{eqnarray}
so by (\ref{eq1}) and (\ref{eq2}) we take that $I=\dfrac{1}{h_0}$ and the proof is completed.

So, we fix some
\begin{eqnarray}
\de_0\in\bigg(0,1-\frac{1}{h_0}\bigg).  \label{eq3}
\end{eqnarray}
After we fix some
\begin{eqnarray}
\el_0\in\bigg(\frac{h_0}{\de_0h_0+1},h_0\bigg),  \label{eq4}
\end{eqnarray}
then $\dfrac{1}{\el_0}<\de_0+\dfrac{1}{h_0}$ and $\el_0>1$.

After, we fix some
\begin{eqnarray}
r_0\in(0,h_0-\el_0).  \label{eq5}
\end{eqnarray}
It holds
\begin{eqnarray}
h_0-\el_0<h_0-1.  \label{eq6}
\end{eqnarray}
After we fix some $N_0\in\N$ such that:
\begin{eqnarray}
(h_0-r_0)^{N_0}>2,  \label{eq7}
\end{eqnarray}
\begin{eqnarray}
\bigg(\frac{h_0-r_0}{\el_0}\bigg)^{N_0}>\frac{2(h_0-\el_0)}{r_0}  \label{eq8}
\end{eqnarray}
\begin{eqnarray}
\el_0^{N_0}>2 \ \ \text{and}  \label{eq9}
\end{eqnarray}
\begin{eqnarray}
\bigg(\frac{8h_0}{h_0-\el_0}\bigg)^{\frac{1}{N_0}}\cdot\frac{1}{\el_0^{\frac{N_0}{N_0+1}}}
<\de_0+\frac{1}{h_0}.  \label{eq10}
\end{eqnarray}
Finally, we fix some positive number $\e_0$ such that:
\begin{eqnarray}
\e_0<\frac{1}{2}\cdot\frac{1}{(2h_0)^{N_0}},  \label{eq11}
\end{eqnarray}
\begin{eqnarray}
\e_0<\frac{r_0}{2},  \label{eq12}
\end{eqnarray}
\begin{eqnarray}
\e_0<\frac{h_0-1}{2}.  \label{eq13}
\end{eqnarray}
We set
\[
K_1:=\oD(h_0+\e_0,\e_0):=\{z\in\C\mid\,\mid z-(h_0+\e_0)\mid\le\e_0\} \ \ \text{and} \ \ L':=K_0\cup K_1.
\]
We prove that for the compact set $L'$ we have $\thi_{L'}<\de_0+\dfrac{1}{h_0}$, where $L'\in C_{h_0}$ of course.

We fix some natural number $n_0>2$.

Let $j_k$, $k=0,1,\ld,n_0$ to be the $n_0$-roots of unity, that is
\[
j_k:=e^{\frac{2k\pi i}{n_0}}, \ \ k=0,1,\ld,n_0-1.
\]
We set
\[
w_k:=(h_0+\e_0)+\e_0j_k, \ \ k=0,1,\ld,n_0-1.
\]
We consider the polynomials
\[
p_1(z):=z^{n_0N_0}-1 \ \ \text{and}
\]
\[
p_2(z):=\prod^{n_0-1}_{k=0}(z-w_k) \ \ \text{and}
\]
\[
p(z):=p_1(z)\cdot p_2(z).
\]
It holds
\begin{eqnarray}
\|p_1\|_{K_0}\le2,  \label{eq14}
\end{eqnarray}
\begin{eqnarray}
\|p_2\|_{K_0}\le(h_0+2\e_0+1)^{n_0}, \ \ \text{so}  \label{eq15}
\end{eqnarray}
\begin{eqnarray}
\|p\|_{K_0}\le2\cdot(h_0+2\e_0+1)^{n_0},  \label{eq16}
\end{eqnarray}
we also take
\begin{eqnarray}
\|p_1\|_{K_1}\le1+(h_0+2\e_0)^{n_0N_0}  \label{eq17}
\end{eqnarray}
\begin{eqnarray}
\|p_2\|_{K_1}\le(2\e_0)^{n_0}, \ \ \text{so} \label{eq18}
\end{eqnarray}
\begin{eqnarray}
\|p\|_{K_1}\le(2\e_0)^{n_0}\cdot(1+(h_0+2\e_0)^{n_0N_0}).  \label{eq19}
\end{eqnarray}
Using the inequalities $r_0<h_0$, (\ref{eq11}) and (\ref{eq12}), we take that
\begin{eqnarray}
2\cdot((2\e_0)\cdot(h_0+2\e_0)^{N_0})^{n_0}\le2\cdot(h_0+2\e_0+1)^{n_0}  \label{eq20}
\end{eqnarray}
and by (\ref{eq16}), (\ref{eq19}) and (\ref{eq20}) we take:
\begin{eqnarray}
\|p\|_{L'}\le2\cdot(h_0+2\e_0+1)^{n_0}.  \label{eq21}
\end{eqnarray}
After we consider the curves $\de_1:\ga_1(t):=\el_0\cdot e^{2\pi it}$, $t\in[0,1]$ and
\[
\de_2:\ga_2(t):=(h_0+\e_0)+r_0e^{2\pi it}, \
 \ t\in[0,1]
\]
we have
\[
|p_1(z)|\ge\el_0^{n_0N_0}-1 \ \ \text{for every} \ \ z\in\de_1
\]
\[
|p_2(z)|\ge(h_0-\el_0)^{n_0} \ \ \text{for} \ \ z\in\de_1.
\]
So, we have
\begin{eqnarray}
|p(z)|\ge(\el_0^{n_0N_0}-1)\cdot(h_0-\el_0)^{n_0}  \ \ \text{for} \ \ z\in\de_1.  \label{eq22}
\end{eqnarray}
We have also
\[
|p_1(z)|\ge(h_0-r_0+\e_0)^{n_0N_0}-1 \ \ \text{for} \ \ z\in\de_2
\]
\[
|p_2(z)|\ge(r_0-\e_0)^{n_0}, \ \ z\in\de_2, \ \ \text{so}
\]
\begin{eqnarray}
|p(z)|\ge((h_0-r_0+\e_0)^{n_0N_0}-1)\cdot(r_0-\e_0)^{n_0} \ \ \text{for} \ \ z\in\de_2.  \label{eq23}
\end{eqnarray}
Using the above inequalities (\ref{eq8}) and (\ref{eq9}) we take that
\begin{eqnarray}
((h_0-r_0+\e_0)^{n_0N_0}-1)\cdot(r_0-\e_0)^{n_0}\ge(\el_0^{n_0N_0}-1)(h_0-\el_0)^{n_0}. \label{eq24}
\end{eqnarray}
We set $\De:=\de_1\cup\de_2$.

By (\ref{eq22}), (\ref{eq23}) and (\ref{eq24}) we take:
\begin{eqnarray}
\min_{z\in\De}|p(z)|\ge(\el_0^{n_0N_0}-1)\cdot(h_0-\el_0)^{n_0}.  \label{eq25}
\end{eqnarray}
Thus, by the inequalities (\ref{eq21}) and (\ref{eq25}) we take that:
\begin{eqnarray}
\frac{\|p\|_{L'}}{\dis\min_{z\in\De}|p(z)|}\le\frac{2\cdot(h_0+2\e_0+1)^{n_0}}
{(\el_0^{n_0N_0}-1)\cdot(h_0-\el_0)^{n_0}}.  \label{eq26}
\end{eqnarray}
By inequalities (\ref{eq9}) and (\ref{eq13}) we take:
\begin{eqnarray}
\frac{2\cdot(h_0+2\e_0+1)^{n_0}}{(\el_0^{n_0N_0}-1)(h_0-\el_0)^{n_0}}<
\frac{2(2h_0)^{n_0}}{\dfrac{1}{2}\el_0^{n_0N_0}\cdot(h_0-\el_0)^{n_0}}. \label{eq27}
\end{eqnarray}
By inequality (\ref{eq10}) we take:
\begin{eqnarray}
\bigg(\frac{2\cdot(2h_0)^{n_0}}{\dfrac{1}{2}\el_0^{n_0N_0}\cdot(h_0-\el_0)^{n_0}}
\bigg)^{\frac{1}{n_0N_0+n_0}}<\de_0+\frac{1}{h_0}.  \label{eq28}
\end{eqnarray}
By inequalities (\ref{eq26}), (\ref{eq27}) and (\ref{eq28}) we have:
\begin{eqnarray}
\bigg(\frac{\|p\|_{L'}}{\min_{z\in\De}|p(z)|}\bigg)^{\frac{1}{n_0N_0+n_0}}<\de_0+\frac{1}{h_0}.
\label{eq29}
\end{eqnarray}
Of course we have $\De\subset\C\setminus L'$. We denote $\OO_{L'}:=(\C\setminus L')\cup\{\infty\}$. It is easy to check that $\De\in\mfd_{L'}$.

We apply now Bernstein's Lemma (5.5.7) of \cite{8} for the polynomial $p$ of degree $n_0N_0+n_0$ and for $z\in\De\subset\OO_{L'}\setminus\{\infty\}$ and we take:
\begin{eqnarray}
\thi_{L',\De}\le\bigg(\frac{\|p\|_{L'}}{\dis\min_{z\in\De}|p(z)|}\bigg)^{\frac{1}{n_0N_0+n_0}}. \label{eq30}
\end{eqnarray}
By (\ref{eq29}) and (\ref{eq30}) we have
\begin{eqnarray}
\thi_{L'}<\de_0+\frac{1}{h},  \label{eq31}
\end{eqnarray}
and the proof of Proposition is complete. \qb
\end{Proof}

We set:
\[
L_0:=K_0\cup\{h_0\}\in C_{h_0}.
\]
We prove now the following proposition.
\begin{prop}\label{prop1.6}
We have
\[
\thi_{L_0}=\frac{1}{h_0}.
\]
\end{prop}
\begin{Proof}
Let some compact set $L_1\in C_{h_0}$. Obviously, we have: $L_0\subset L_1$.

We set: $\OO_0:=(\C\setminus L_0)\cup\{\infty\}$ and $\OO_1:=(\C\setminus L_1)\cup\{\infty\}$. Obviously, $\OO_1\subset\OO_0$. By Corollary 4.4.5 \cite{8} we have:
\[
g_{\OO_1}(z,\infty)\le g_{\OO_0}(z,\infty), \ \ z\in\OO_1.
\]
So
\setcounter{equation}{0}
\begin{eqnarray}
e^{-g_{\OO_0}(z,\infty)}\le e^{-g_{\OO_1}(z,\infty)}\Rightarrow
\thi_{L_0,\De}\le\thi_{L_1,\De}, \ \ \text{for every} \ \ \De\in\mfd_{L_1}.  \label{eq1}
\end{eqnarray}
Of course we have $\mfd_{L_1}\subset\mfd_{L_0}$. We have
\begin{align}
\{\thi_{L_0,\De\mid\De\in\mfd_{L_1}}\}\subseteq\{\thi_{L_0,\De,\mid\De\in\mfd_{L_0}}\}
&\Rightarrow\inf\{\thi_{L_0,\De\mid\De\in\mfd_{L_0}}\}\le\inf\{\thi_{L_0,\De\mid
\De\in\mfd_1}\}\nonumber\\
&\Rightarrow\thi_{L_0}\le\inf\{\thi_{L_0,\De\mid\De\in\mfd_{L_1}}\}.  \label{eq2}
\end{align}
By (\ref{eq1}) we have
\begin{eqnarray}
\inf\{\thi_{L_0,\De}\mid\De\in\mfd_{L_1}\}\le\inf\{\thi_{L_1,\De}\mid\De\in\mfd_1\}=\thi_{L_1}. \label{eq3}
\end{eqnarray}
By (\ref{eq2}) and (\ref{eq3}) we have:
\begin{eqnarray}
\thi_{L_0}\le\thi_{L_1}.  \label{eq4}
\end{eqnarray}
By (\ref{eq4}) and Proposition \ref{prop1.5} we have
\begin{eqnarray}
\thi_{L_0}\le\frac{1}{h_0}.  \label{eq5}
\end{eqnarray}
By the proof of Proposition \ref{prop1.5}, it is easy to see that we can construct a strictly decreasing sequence of compact sets $L_n\in C_{h_0}$, such that $n>\dfrac{1}{1-\dfrac{1}{h_0}}$, $L_{n+1}\subset L_n$, for $n>\dfrac{1}{1-\dfrac{1}{h_0}}$, and $\thi_{L_n}<\dfrac{1}{n}+\dfrac{1}{h_0}$. Of course $\bigcap_{n>a}L_n=L_0$, where $a:=\dfrac{1}{1-\dfrac{1}{h_0}}$. We have $\Big(\bigcap\limits_{n>a}L_n\Big)^c=L^c_0\Rightarrow\bigcup\limits_{n>a}L^c_n=L^c_0$. We set
\[
\OO_n:=(\C\setminus L_n)\cup\{\infty\}, \ \ n>a,
\]
so, we have $\OO_0=\bigcup\limits_{n>a}\OO_n$.

We fix some $\De_0\in\mfd_{L_0}$. It is easy to see that there exists some $m_0>a$ such that $\De_0\in\mfd_{L_n}$, for every $n\ge m_0$. Of course we have $\OO_0=\bigcup\limits_{n\ge m_0}\OO$.

By Theorem 4.4.6 \cite{8} we have
\[
\lim_{n\ra\infty}g_{\OO_n}(z,\infty)=g_{\OO_0}(z,\infty), \ \ z\in\OO_0, \ \ \text{so}
\]
\[
\lim_{n\ra\infty}g_{\OO_n}(z,\infty)=g_{\OO_0}(z,\infty), \ \ \text{for every} \ \ z\in\De_0.
\]
Of course
\[
\OO_n\subset\OO_{n+1}, \ \ n\ge m_0, \ \ \text{so}
\]
\[
g_{\OO_n}(z,\infty)\le g_{\OO_{n+1}}(z,\infty), \ \ n\ge m_0.
\]
This gives that the sequence of functions $\{-g_{\OO_n},\;n\ge n_0\}$ is a decreasing sequence of continuous functions on the compact set $\De_0$, so by Dini' Theorem we have that $g_{\OO_n}\ra g_{\OO_0}$ on $\De_0$ uniformly.

We fix some positive number $\e_0$. Then there exists some $m_1\ge m_0$ such that
\[
|g_{\OO_n}(z,\infty)-g_{\OO_0}(z,\infty)|<\e_0 \ \ \text{for every} \ \ z\in\De_0, \ \ n\ge m_1,
\]
so we take
\begin{align*}
&g_{\OO_0}(z,\infty)-g_{\OO_n}(z,\infty)<\e_0, \ \ z\in\De_0, \ \ n\ge m_1 \\
\Rightarrow\,&e^{g_{\OO_0}(z,\infty)}<e^{\e_0}\cdot e^{g_{\OO_n}(z,\infty)}, \ \ z\in\De_0, \ \ n\ge m_1 \\
\Rightarrow\,&e^{-\e_0-g_{\OO_n}(z,\infty)}<e^{-g_{\OO_0}(z,\infty)}, \ \ z\in\De, \ \ n\ge m_1 \\
\Rightarrow\,&e^{-\e_0}\thi_{L_n,\De_0}<\thi_{L_0,\De_0}, \ \ n\ge m_1 \\
\Rightarrow\,&e^{-\e_0}\thi_{L_n}<\thi_{L_0,\De_0}, \ \ n\ge m_1 \\
\Rightarrow\,& e^{-\e_0}\frac{1}{h_0}<\thi_{L_0,\De_0}.
\end{align*}
This holds for every $\e>0$, so $\dfrac{1}{h_0}\le\thi_{L_0,\De_0}$. This inequality holds for every $\De\in\mfd_{L_0}$, so
\begin{eqnarray}
\frac{1}{h_0}\le\thi_{L_0}.  \label{eq6}
\end{eqnarray}
By (\ref{eq5}) and (\ref{eq6}) we have $\thi_{L_0}=\dfrac{1}{h_0}$ and the proposition is complete. \qb
\end{Proof}
\section{Final step of the proof of Theorem \ref{mainthm}: $\bbb{\thi_L}=\bbb{\rho_L}$}\label{sec2}
\noindent

So, by the above Theorem \ref{thm1.4} we have proved that the number $\thi_L$ is positive and we have found an easy-computed (in all simple cases) lower bound of $\thi_L$.

For the sequel, we refer to \cite{8} for the respective terminology.

More specifically:

For the definition of Harnack distance see Definition 1.3.4. We note that the Harnack distance is a continuous function. For the definition of logarithmic capacity see Definition 5.1.1. For the definition of a Fekete $n$-tuple and the $n$-th diameter $\de_n(K)$ of a compact set $K\subseteq\C$ see Definition 5.5.1.

For the definition of a Fekete polynomial of degree $n\ge2$ see Definition 5.5.3, for some compact set $K$. If $A\subseteq\C$ and $F:A\ra\C$ be a complex function we denote
\[
\|F\|_A:=\sup\{x\in\R\mid\exists\,a\in A:x=|F(a)|\}\in[0,+\infty].
\]

We remind here (Bernstein's Lemma) Theorem 5.5.7 of \cite{8}.

Let $L$ be a non-polar compact subset of $\C$, and let $\OO$ be the component of\\ $(\C\cup\{\infty\})\sm L$ containing $\infty$. If $q_n$ is a Fekete polynomial of degree $n\ge2$, for $L$ then
\[
\bigg(\frac{|q_n(z)|}{\|q_n\|_L}\bigg)^{1/n}\ge e^{g_\OO(z,\infty)}\bigg(\frac{c(L)}{\de_n(L)}\bigg)^{T_\OO(z,\infty)} \ \ \text{for every} \ \ z\in\OO\sm\{\infty\}  \eqno{\mbox{($\ast$)}}
\]
We consider now the fixed set $L$ of our work, where $L=\bigcup\limits^{m_0}_{i=0}K_i$, $m_0>1$, $\overset{\circ}{K}_0\neq\emptyset$ and $K_i$, $i=1,\ld,m_0$, contains more than one point, $K_i$, $i=0,1,\ld,m_0$ are the connected components of $L$ and $L^c$ is connected.

We choose some Fekete polynomial $q_m$, for every $m=2,3,\ld$, for $L$, and we fix them for the sequel.

We set
\[
\inf_{z\in\De}|q_m(z)|:=\inf\{x\in\R\mid\exists\,z\in\De:x=|q_m(z)|\} \ \ \text{for every} \ \ m=2,3,\ld, \ \ \De\in\mfd_L.
\]
By the above terminology we get the following lemma using inequality $(\ast)$.
\begin{lem}\label{lem2.1}
For every positive constant $c\in(\thi_L,1)$ there exists some $\De\in\mfd_L$ and some natural number $m_{L,\De,c}$ that depends  only on $L,\De$ and $c$ such that:
\[
\frac{\|q_n\|_L}{\dis\inf_{z\in\De}|q_n(z)|}<c^n, \ \ \text{for every} \ \ n\ge m_{L,\De,c}.
\]
\end{lem}
\begin{Proof}
Take arbitrary $c\in(\thi_L,1)$. By the definition of the number $\thi_L$ we can take some $\De\in\mfd_L$ such that $\thi_L<\thi_{L,\De}<c$, where $\De$ depends on $L,c$.

We fix some $m_0\ge2$. By $(\ast)$ we get
\[
\bigg(\frac{\|q_{m_0}\|_L}{|q_{m_0}(z)|}\bigg)^{1/m_0}\le\frac{1}{e^{g_\OO(z,\infty)}}
\bigg(\frac{\de_{m_0}(L)}{c(L)}\bigg)^{T_\OO(z,\infty)} \ \ \text{for every} \ \ z\in\De.
\]
So we have
\setcounter{equation}{0}
\begin{eqnarray}
\bigg(\frac{\|q_{m_0}\|_L}{|q_{m_0}(z)|}\bigg)^{1/m_0}\le\thi_{L,\De}\cdot
\bigg(\frac{\de_{m_0}(L)}{c(L)}\bigg)^{\|T_\OO\|_\De} \label{eq1}
\end{eqnarray}
By (\ref{eq1}) we have:
\begin{eqnarray}
\bigg(\frac{\|q_{m_0}\|_L}{\dis\inf_{z\in\De}|q_{m_0}(z)|}\bigg)^{1/m_0}\le\thi_{L,\De}\cdot
\bigg(\frac{\de_{m_0}(L)}{c(L)}\bigg)^{\|T_\OO\|_\De}  \ \ \text{for every} \ \ m_0\ge2.  \label{eq2}
\end{eqnarray}
By Fekete-Szeg\"{o} Theorem (Theorem 5.5.2) we have that
\begin{eqnarray}
\de_m(L)\ra c(L) \ \ \text{as} \ \ m\ra+\infty.  \label{eq3}
\end{eqnarray}
The number $\thi_{L,\De}$ depends only on $\De$ and $L$.

Thus, because $\thi_{L,\De}<c$. (4) there exists some natural number $m(0,\De)=m_\De$ that depends on $\De$ such that
\[
\frac{\|q_m\|_L}{\dis\inf_{z\in\De}|q_m(z)|}<c^m \ \ \text{for every} \ \ m\in\N, \ \ m\ge m_\De,
\]
by (\ref{eq2}), (\ref{eq3}) and (\ref{eq4}).

This completes the proof. \qb
\end{Proof}

We will need also a proposition that is a variation of the well known Bernstein-Walsh Theorem.
\begin{prop}\label{prop2.2}
Let some compact set $L=\bigcup\limits^{m_0}_{i=0}K_i$, $m_0\in\N$, as above where $K_i$, $i=0,1,\ld,m_0$ are the connected components of $L$. Let some complex polynomials $p_j$, $j=0,1,\ld,m_0$. We consider the function $F:L\ra\C$ that is defined by the following formula:
\[
F(z)=p_j(z) \ \ \text{if} \ \ z\in K_j \ \ \text{for every} \ \ j=0,1,\ld,m_0.
\]
Then, for every positive number $c\in(\thi_L,1)$, there exists $\De\in\mfd_L$, that depends on $L,c$, some natural number $m=m_\De$ that depends on $\De$, some positive constant $A=A_{\De,F}$ that depends on $\De,F$ and some sequence of polynomials $(S_j)$ that depends on $\De,F$ so that the following inequality holds:
\[
\|F-S_m\|_L<A\cdot c^m \ \ \text{for every} \ \ m\in\N, \ \ m\ge m_{\De}, \ \ \deg(S_m)\le m-1.
\]
\end{prop}
\begin{Proof}
We consider the compact set $L=\bigcup\limits^{m_0}_{i=0}K_i$, the polynomials $p_j$, $j=0,1,\ld,m_0$ and the complex function $F$ as in the suppositions of this proposition. We fix some positive number $c_0\in(\thi_L,1)$. After we fix some $\De=\De_{L,c_0}\in\mfd_L$ that depends on $L,c_0$ such that $\thi_{L,\De}\in(\thi_L,c_0)$.

Afterwards we apply Lemma \ref{lem2.1} and we get that there exists some natural number $m_\De>2$ that depends on $\De$ such that:
\setcounter{equation}{0}
\begin{eqnarray}
\frac{\|q_n\|_L}{\dis\inf_{z\in\De}|q_n(z)|}<c^n_0, \ \ \text{for every} \ \ n\in\N \ \ n\ge m_\De>2.  \label{eq1}
\end{eqnarray}
Let $\De=\bigcup\limits^{m_0}_{i=0}\de_i$ where $\de_i$, $i=0,1,\ld,m_0$ are the connected components of $\De$. There exist $V_i$, $i=0,1,\ld,m_0$ bounded simply connected domains pairwise disjoint such that $\de_i\subset V_i$ for every $i=0,1,\ld,m_0$ by the proof of Lemma \ref{lem1.2}.

We set $V:=\bigcup\limits^{m_0}_{i=0}V_i$.

We consider the function $\Fi:V\ra\C$ by the following formula:
\[
\Fi(z)=p_j(z) \ \ \text{for every} \ \ z\in V_j \ \ \text{for every} \ \ j=0,1,\ld,m_0.
\]
Of course we have $\Fi\upharpoonright_L=F$, and $\Fi$ is holomorphic.

We consider some fixed sequence $(q_m)$, $m\ge2$ of Fekete polynomials for $L$ of degree $m$, $m\ge2$. We define now the functions $r_m:L\ra\C$ with the formula:
\begin{eqnarray}
r_m(w):=\sum^{m_0}_{j=0}\frac{1}{2\pi i}\int\limits_{\de_j}\frac{\Fi(z)}{q_m(z)}\cdot
\frac{q_m(w)-q_m(z)}{w-z}dz  \label{eq2}
\end{eqnarray}
for every $m\in\N$, $m\ge2$, $w\in L$. For every $m\ge2$, the functions $r_m$ are polynomials of degree at most $m-1$ as we can see easily.

It is obvious that $\ssum^{m_0}_{j=0}\tind_{\de_j}(a)=0$ for every $a\in\C\sm V$.

We fix some $n_0\ge2$.

We apply now the global Cauchy's integral formula for the function $\Fi$ and the smooth Jordan curves $\de_j$, $j=0,1,\ld,m_0$ and we take:
\begin{eqnarray}
\sum^{m_0}_{j=0}\tind_{\de_j}(w)\cdot\Fi(w)=\sum^{m_0}_{j=0}\frac{1}{2\pi i}\int\limits_{\de_j}
\frac{\Fi(z)}{z-w}dz \ \ \text{for every} \ \ w\in V\sm\De.  \label{eq3}
\end{eqnarray}
By (\ref{eq2}) and (\ref{eq3}) we get:
\begin{eqnarray}
\Fi(w)-r_{n_0}(w)=\sum^{m_0}_{j=0}\frac{1}{2\pi i}\int\limits_{\de_j}\frac{\Fi(z)q_{n_0}(w)}{(z-w)q_{n_0}(z)}dz \ \ \text{for every} \ \ w\in L.  \label{eq4}
\end{eqnarray}
We set $\el_j=\text{length}(\de_j)$ for every $j=0,1,\ld,m_0$, $\la_0:=\ssum^{m_0}_{j=0}\el_j$,
\[
dist(\De,L):=\min\{x\in\R\mid\exists\,z_1\in\De,z_2\in L:x=|z_1-z_2|\}.
\]
After we define the number:
\begin{eqnarray}
A:=\frac{\la_0\cdot\|\Fi\|_L}{2\pi\cdot dist(\De,L)}.  \label{eq5}
\end{eqnarray}
Of course the above number $A$ depends on $F,\De$. By (\ref{eq4}) and (\ref{eq5}) we take easily that:
\begin{eqnarray}
\|F-r_{n_0}\|_L\le A\cdot\frac{\|q_{n_0}\|_L}{\dis\inf_{z\in\De}|q_{n_0}(z)|}.  \label{sex26}
\end{eqnarray}
Of course, the positive number $A$ is independent from the natural number $n_0$.

So by (\ref{sex26}) we get:
\begin{eqnarray}
\|F-r_n\|_L\le A\cdot\frac{\|q_n\|_L}{\dis\inf_{z\in\De}|q_n(z)|} \ \ \text{for every} \ \ n\in\N, \ \ n\ge2.  \label{sex27}
\end{eqnarray}
By (\ref{eq1}) and (\ref{sex27}) we get that:
\begin{eqnarray}
\|F-r_n\|_L<A\cdot c^n_0 \ \ \text{for every} \ \ n\in\N, \ \ n\ge m_{L,\De,c_0},  \label{sex28}
\end{eqnarray}
where the natural number $m_\De$ depends on $\De$, the set $\De$ depends on $L,c_0$, the positive number $c_0$ depends on $L$, the constant $A$ depends on $F,\De$ and the polynomials $r_n,n\ge2$ depend on $F,\De$.

The above inequality (\ref{sex28}) completes the proof of this proposition for $S_m=r_m$, $m\ge2$. \qb
\end{Proof}

The above Proposition \ref{prop2.2} gives some important role to the number $\thi_L$. It shows that the number $\thi_L$ plays a crucial role in the problem of approximation by polynomials.

We connect here the number $\thi_L$ with an other number that is a characteristic for the compact set $L$. This number is the well known number that is called the asymptotic convergence factor and is noted $\rho_L$ as usually, see \cite{6}, \cite{10} and \cite{11}.

First of all we define here the number $\rho_L$.

Let $L:=\bigcup\limits^{m_0}_{i=0}K_i$, $m_0\in\N$ be some compact set and $K_i$, $i=0,1,\ld,m_0$, are simple compact sets pairwise disjoint.
Let
\[
A(L):=\{f:L\ra\C\mid f \ \ \text{is continuous on} \ \ L \ \ \text{and holomorphic on} \ \ \lo\}.
\]
We consider the space $A(L)$ endowed with the supremum norm $\|\cdot\|_\infty$. As it is well known the space $(A(L),\|\cdot\|_\infty)$ is a Banach Algebra. We fix some polynomials $p_j$, $j=0,1,\ld,m_0$ such that
$p_i\neq p_j$ for every $i,j\in\{0,1,\ld,m_0\}$, $i\neq j$.

We consider the function $f:L\ra\C$, that is defined by the following formula:
\[
f(z)=p_j(z) \ \ \text{if} \ \ z\in K_j \ \ \text{for every} \ \ j=0,1,\ld,m_0.
\]
Of course $f\in A(L)$ but $f$ is not a polynomial.

Let $n\in\N$, $n\ge1$ be some natural number and
\[
V_n:=\{q:L\ra\C\mid q  \ \  \text{is a polynomial of degree at most}\ \  n\}.
\]

We note $t_{f,n}:=dist(V_n,f)$ for $n=1,2,\ld\;.$ Of course $V_n\subset A(L)$ for every $n=1,2,\ld,\ld$ and the space $V_n$ is a closed vector proper subspace of $A(L)$.

(Of course the space $A(L)$ has a non-denumerable hamel basis whereas the space $V_n$ has dimension $n$ for every $n=1,2,\ld\;.$). The number $t_{f,n}$ is positive for every $n=1,2,\ld$ as the distance of the closed subset $V_n$ of $A(L)$ from the compact set $\{f\}\subset A(L)$ where $f\notin V_n$ for every $n=1,2,\ld\;.$ It is well known that there exists some $h_n\in V_n$ such that $\|f-h_n\|=t_{f,n}$ for every $n=1,2,\ld$ and this polynomial $h_n$ is unique. (See \cite{14}).

Of course, the sequence $t_{f,n}$ is decreasing and by Mergelyan's Theorem (or Runge's Theorem) \cite{9} we have that $\dis\lim_{n\ra+\infty}t_{f,n}=0$.

But Mergelyan's Theorem (or Runge's Theorem) does not tell us something about the rate of convergence of the sequence $t_{f,n}$ is general.

However, in this specific case where $f$ is a polynomial on every compact set $K_i$, $i=0,1,\ld,m_0$, we have by Bernstein-Walsh Theorem and our Proposition \ref{prop2.2} information about the rate of convergence of the sequence $t_{f,n}$ that tends to zero with an exponential rate, that is very fast.

Now, it is well known that there exists some positive number $\rho_L\in(0,1)$ such that: $\underset{n\ra+\infty}{\lim\sup}t^{1/n}_{f,n}=\rho_L$ (see \cite{14}, \cite{6}).

The number $\rho_L$ depends only on the compact set $L$ and is independent from the specific function $f$.

The number $\rho_L$ is called the asymptotic convergence factor for the compact set $L$.

We have the following very important information about the number $\rho_L$ below.
\begin{prop}\label{prop2.3}
By the previous notations we have that: $\rho_L=\thi_L$.
\end{prop}
\begin{Proof}
We take some $\De\in\mfd_L$.

Let $\De=\bigcup\limits^{m_0}_{i=0}\de_i$ where $\de_i$, $i=0,1,\ld,m_0$ be the connected components of $\De$.

Let $G_i$, $i=0,1,\ld,m_0$ be bounded simply connected domains, pairwise disjoint such that $\de_i\subset G_i$ for $i=0,1,\ld,m_0$, (using the proof of Lemma \ref{lem1.2}).

We consider $m_0+1$ polynomials $p_i$, $i=0,1,\ld,m_0$ where $p_i\neq p_j$ for every $i,j\in\{0,1,\ld,m_0\}$, $i\neq j$. We define the holomorphic function $F:G\ra\C$, where $G:=\bigcup\limits^{m_0}_{i=0}G_i$ with the following formula:
\[
F(z)=p_i(z) \ \ \text{for every} \ \ z\in G_i, \ \ i=0,1,\ld,m_0.
\]
We apply Proposition \ref{prop2.2} and for the above function $F$ there exists a sequence of polynomials $r_n$, $n\ge2$, some positive number $A$ and some natural number $n_0$ such that
\[
\|F-r_n\|_L<A\cdot c^n_0 \ \ \text{for every} \ \ n\ge n_0
\]
for some positive constant $c_0\in(\thi_{L,\De},1)$, (see the proof of Lemma \ref{lem2.1} also).

This gives that:
\setcounter{equation}{0}
\begin{eqnarray}
\underset{n\ra+\infty}{\lim\sup}\|F-r_n\|^{1/n}_L\le c_0.  \label{eq1}
\end{eqnarray}
Let $S_n$, $n=2,3,\ld$ be the unique polynomial of degree at most $n$ (that there exists see \cite{14}) that minimizes the quantity $\|F-S_n\|_L$. We write $t_{F,n}:=\|F-S_n\|$ for simplicity. It is known \cite{14}, \cite{6},
\begin{eqnarray}
\underset{n\ra+\infty}{\lim\sup}t^{1/n}_{F,n}=\rho_L.  \label{eq2}
\end{eqnarray}
By the definition of the number $t_{F,n}$ we have of course:
\begin{eqnarray}
t_{F,n}\le\|F-r_n\| \ \ \text{for} \ \ n\ge2.  \label{eq3}
\end{eqnarray}
By (\ref{eq1}), (\ref{eq2}) and (\ref{eq3}) we get: $\rho_L\le c_0$. But the number $c_0$ is some arbitrary positive number such that $\thi_{L,\De}<c_0<1$. This gives that $\rho_L\le\thi_{L,\De}$. Because this holds for every $\De\in\mfd_L$ we get
\begin{eqnarray}
\rho_L\le\thi_L.  \label{eq4}
\end{eqnarray}
Now $\rho_L=\exp(-g_c)$ where $g_c$ is the critical potential (see \cite{6}) and\\ $\ga:=\{z\in\C:g_\OO(z)=g_c\}$ is the critical level curve where $\OO:=(\C\cup\{\infty\})\sm L$ and $g_\OO$ is the Green's function for $L$. It is simple to see by the continuity of $g_\OO$ that there exists a sequence of curves $\De_n$, $n=1,2,\ld$, where $\De_n\in\mfd_L$ for $n=1,2,\ld$ such that
\begin{eqnarray}
\thi_{L,\De_n}\ra\rho_L \ \ \text{as} \ \  n\ra+\infty.  \label{eq5}
\end{eqnarray}
By (\ref{eq4}) and (\ref{eq5}) we obtain that $\rho_L=\thi_L$ and the proof of this proposition is\linebreak complete. \qb
\end{Proof}

Replacing in Theorem \ref{thm1.4} the compact set $K_0$ by any $K_j$ for $j=1,\ldots ,m$ and using Proposition \ref{prop2.3} the proof of Theorem \ref{mainthm} is complete.

%
%
%


%
%

Department of Mathematics and Applied Mathematics, University of Crete, Panepistimiopolis Voutes, 700-13, Heraklion, Crete, Greece.\\
email:tsirivas@uoc.gr


\begin{thebibliography}{99}
%
\bibitem{1} N. Akhiezer, Theory of Approximation, {\em Ungar}, New York, (1956).
%
\bibitem{2}  D. H. Armitage. S. J. Gardiner, Classical Potential Theory, {\em Springer}, London, (2001).
    %
\bibitem{3} R. B. Burckel, An Introduction to Classical Complex Analysis, {\em Birkh\"{a}user Verlag}, Basel, (1979).
%
\bibitem{4} E. W. Cheney, Introduction to Approximation Theory, {\em McGraw-Hill}, New York, (1966).
%
\bibitem{5} G. Costakis, Some remarks on universal functions and Taylor series, {\em Math. Proc. Cambridge-Philos. Soc.} {\bf128} (2000), 157-175.
%
\bibitem{6} M. Embree, L. N. Trefethen, Green's functions for multiply connected domains via conformal mapping. {\em SIAM Rev.}, {\bf41}, 745-761, (1999).
%
\bibitem{7} K.-G. Grosse-Erdmann, Holomorphe Monster und universelle Functionen, {\em Mitt. Math. Sem. Giessen}, {\bf176} (1987).
%
\bibitem{8} T. Ransford, Potential Theory in the Complex Plane, {\em Cambridge Univ. Press}, Cambridge, (1995).
%
\bibitem{9}  W. Rudin, Real and Complex Analysis, 3rd ed. {\em McGraw-Hill}, (1966).
%
\bibitem{10} Klaus Schiefermayr, Estimates for the Asymptotic Convergence Factor of Two Intervals, {\em Journal of Computational and Applied Mathematics}, {\bf236}, 26-36, (2011).
%
\bibitem{11} Tobin A. Driscoll, Kim-chuan Toh, Lloyd Trefethen, From potential theory to matrix iterations in six steps, {\em Journal SIAM} Review, volume {\bf40}, 547-578, (1998).
%
%
\bibitem{13} N. Tsirivas, Universal Taylor series on specific compact sets, submitted.
%
\bibitem{14} J. L. Walsh. Interpolation and Approximation by Rational functions in the Complex domain, 5th ed., {\em Amer. Math. Soc.} Providence, RI, (1969).

\end{thebibliography}
\end{document}